\documentstyle[12pt]{amsart}
\theoremstyle{plain}
\newtheorem{Thm}{Theorem}[section]
\newtheorem{Cor}[Thm]{Corollary}
\newtheorem{Prop}[Thm]{Proposition}
\newtheorem{Lem}[Thm]{Lemma}

\newtheorem{Def}[Thm]{Definition}
\theoremstyle{remark}

\newtheorem{Conj}[Thm]{Conjecture}
\newtheorem{Rec}[Thm]{Recall}
\title{Green's Conjecture for the generic canonical curve.}
\author{ Montserrat\ Teixidor i Bigas}
\address
{Mathematics Department\\
Tufts University\\ Medford MA 02155\\ U.S.A.}
\email{teixidor@@dpmms.cam.ac.uk\\mteixido@@tufts.edu}

\begin{document}
\maketitle

\section*{Introduction}
Let $C$ be a non-singular curve of genus $g$ over an algebraically
closed field ${\bf k}$ of characteristic zero. Let $K$ be the canonical 
sheaf on $C$. If $C$ is not hyperelliptic, the map associated to the
 complete canonical series $|K|$
$$C\rightarrow {\bf P}^{g-1}$$
is an embedding and the image curve is projectively normal. If the
 curve is neither trigonal nor a plane quintic, the ideal of $C$
is known to be generated by quadrics. Continuing in this vein,
Mark Green made the conjecture that the resolution of the ideal
of $C$ in ${\bf P}^{g-1}$ should depend on the linear series that $C$ 
has.

To make this precise, one defines property $N_p$.
Take a minimal resolution of the ideal sheaf of $C$ in ${\bf P}^{g-1}$.
Then one says that  property $N_0$ holds if $C$ is projectively
normal, $N_1$ means that the ideal of the curve is generated 
by quadrics, $N_2$ means that in addition the syzygies among these quadrics
are generated by linear relations...In general $N_p$ means that $N_{p-1}$
holds and the $p^{th}$ syzygies are generated by linear relations.

Define the Clifford index of $C$ by 
$$Cliff(C)=min\{ degL-2(h^0(C,L)-1)|L\in Pic C, h^0(C,L)\ge 2,h^1(C,L)\ge 2\}$$
In particular, for a curve that is generic in the sense of moduli, the 
Clifford index of $C$ is given by $[(g-1)/2]$ while the most special curves
from the point of view of Clifford index are hyperelliptic curves that have
$Cliff(C)=0$.
Then 
\begin{Conj}
\label{conj}
[Green](cf. \cite{G} 5.1)
 The curve $C$ has property $N_p$ if and only if 
$Cliff (C)>p$.
\end{Conj}

The only if part of \ref{conj} was proved by Green and Lazarsfeld
 (cf. \cite{G} Appendix),
The conjecture has been proved for $g\le 8$ (cf \cite{S1}) and 
for $p=2$ (i.e. $C$ quatrigonal) (cf \cite{V,S2,HPR}).

In a slightly more modest vein
\begin{Conj} 
\label{conjgen}[Generic Green's conjecture]
(cf. \cite{G} 5.6)
 The generic curve $C$ of genus $g$ satisfies $N_{[(g-3)/2]}$.
\end{Conj}

Few effective results seem to be available with respect to \ref{conjgen}.
The only published work that we are aware of is \cite{E} 
where L.Ein  showed that the first two steps of
the resolution of a generic canonical curve are as expected.
The relevance of \ref{conjgen} is enhaced by the following result of 
Hirschowitz and Ramanan (cf. \cite{HR} Theorem 1.1).

\begin{Thm}[Hirschowitz-Ramanan]
\label{HR}
For odd $g=2k+1\ge 5$, Green's Conjecture holds for the generic 
curve if and only if it holds for all curves of (maximal) 
Clifford index $k$.
\end{Thm}

The purpose of this paper is to prove the generic Green's Conjecture.

\begin{Thm}
\label{teorema}
Let $C$ be a generic curve of genus $g$. Then, Green's conjecture on
 the syzygies of  the canonical curve $C$ holds.
\end{Thm}

As a consequence of Hirschowitz and Ramanan's Theorem,
 one obtains a very concrete open dense subset of ${\cal M}_g$
for odd $g$ where Green's Conjecture actually holds. In the language of 
Hirschowitz and Ramanan, this is the specific Green's Conjecture

\begin{Cor}
If $g=2k+1\ge 5$, then Green's Conjecture holds for all curves of 
maximal Clifford index $k$ (i.e. $C$ satisfies $N_{k-1}$ if and only if 
$Cliff(C)=k$).
\end{Cor}

The proof of \ref{teorema} is as follows: Define a vector bundle $E$
 as the dual of the kernel of the evaluation map of the canonical 
linear series. Namely, $E^*$ is defined by the exact
sequence
\begin{Def}
\label{E}
$$0\rightarrow E^*\rightarrow (H^0(C,K_C))^*
\otimes {\cal O}_C\rightarrow K \rightarrow 0.$$
\end{Def}

From a result of Paranjape-Ramanan (cf. \cite{PR} Remark 2.8, p.507), 
Green's Conjecture would follow from the surjectivity of the maps
$$\wedge ^r(H^0(C,K_C))^* \rightarrow H^0(\wedge ^r E), r\le Cliff(C).$$
 We want to prove that
 this is the case for $C$ generic. To this end consider a hyperelliptic
curve $C_0$. Notice that the map  
$$(H^0(C,K_{C_0}))^*\otimes {\cal O}_{C_0}\rightarrow E \rightarrow 0$$
identifies $(H^0(C,K_{C_0}))^*  $ to a subspace $W$ of $H^0(C_0,E)$.
 As $C_0$ is
 hyperelliptic, this subspace is proper. We shall start by computing 
$W$ and the image $W^r\subset H^0(C_0,\wedge ^rE)$ of its exterior
 powers $\wedge ^r W$. Every infinitessimal
 defformation of the curve $C_0$
 determines a unique  infinitessimal defformation of $E$ 
that preserves it as the dual of the kernel of the canonical evaluation 
map. Consider the defformation of $\wedge ^r E$ that this induces.
 We then see that the only sections of 
$H^0(\wedge ^rE)$ that give rise to sections of the infinitessimal 
defformation of $\wedge ^rE$ are those in $W^r$. This will conclude
the proof.

Acnowledgments: I would like to thank the following people and institutions
that contributed in different ways to the present work. 
My interest on Green's Conjecture came from a conversation with S.Ramanan.
 Mark
Green read a first version of the paper and made several 
suggestions to improve the presentation.
I am visiting the Pure Mathematics Department at the University of 
Cambridge and I benefitted from conversations with Tony Wasserman and 
Nick Shepherd-Barrow as well as e-mails with Loring Tu.
 I am a member of the Europroj group
Vector Bundles on Algebraic Curves. 

\section{Identification of the vector bundle $E$ and its space of sections}
{\bf Notations} In this section and section 3, $C=C_0$ denotes a hyperelliptic
 curve  and $L$ the unique line bundle of degree
 two with two sections.

If $F$ is a vector bundle on a curve $C$, we write $H^i(F)$ for $H^i(C,F)$
if there is no danger of confusion.
 If several curves are involved or if we are
considering sections on an open set only, we shall make this clear.

As in the previous section $E$ will denote the vector bundle defined 
in \ref{E}, $W=W^1$ will be the image of $H^0(K)^*\subset H^0(E)$ and 
$W^r$ the image of the exterior powers of $W$ in $H^0(\wedge^r E)$.

\begin{Prop}
\label{W}
 Let $\pi$ be the map $C_0\rightarrow {\bf P}^1$ associated to 
$L$. Let $H$ be the space of sections of $L$. Then, 
$$E\cong S^{g-2}H^*\otimes \wedge ^2H^*\otimes L$$
Moreover, the map $(H^0(K))^*\rightarrow H^0(E)$ can be identified to the 
natural inclusion
$$\varphi _1 :S^{g-1}H^*\rightarrow S^{g-2}H^*\otimes \wedge ^2H^*\otimes H$$
\end{Prop}
\begin{pf}
Note that the canonical sheaf $K_{C}$ on $C$ is of the form 
$$K=L^{\otimes g-1}=\pi ^*({\cal O}_{\bf P^1}(g-1))$$ 
and
 $$H^0(C,K_C)=S^{g-1}(H)
=\pi ^*(H^0({\bf P}^1,{\cal O}_{\bf P^1}(g-1)))$$
Therefore, the exact sequence defining $E^*$ is the pull-back of
 the exact sequence in ${\bf P}^1$
\begin{equation}
\tag{*}0\rightarrow Ker \rightarrow 
H^0({\bf P}^1,{\cal O}_{\bf P^1}(g-1))\otimes {\cal O}_{{\bf P}^1}
\rightarrow {\cal O}_{\bf P^1}(g-1))\rightarrow 0 \end{equation}

Denote by $\bar H$ the space of sections $H^0({\bf P^1},
 {\cal O}_{{\bf P^1}}(1))$. 
Tensoring (*) with ${\cal O}(1)$ and taking global sections,
one obtains
$$0\rightarrow H^0(Ker\otimes {\cal O}(1)) \rightarrow S^{g-1}\bar H\otimes
 \bar H \rightarrow S^g \bar H\rightarrow 0$$
Hence, $H^0(Ker\otimes {\cal O}(1))\equiv S^{g-2}\bar H\otimes
 \wedge ^2\bar H$ (cf. for example \cite{FH}  p.224 15.20). Using 
the exact sequence (*), one checks that $Ker$ is a direct sum of line bundles 
of degree $-1$. Hence $Ker= S^{g-2}\bar H\otimes \wedge ^2\bar H\otimes
{\cal O}(-1)$. As $E^*=\pi^*(Ker), H=\pi ^*\bar H$, the result follows.

\end{pf}

\begin{Prop}
\label{W^r}
Let the notations be as above. Choose $r\le g-1$
Then,
 $$\wedge ^rE=\wedge ^r(S^{g-2}H^*\otimes \wedge ^2H^*)\otimes L^{\otimes r}$$
Moreover, the map
$$\wedge^rH^0(E) \rightarrow H^0(\wedge ^r E)$$
can be identified to 
$$\varphi _r :  \wedge^r (S^{g-1}H^*)\rightarrow \wedge ^r(S^{g-2}H^*
\otimes \wedge ^2H^*)\otimes S^rH$$
and is an immersion.
 In particular, its image $W^r$ has dimension
 $\begin{pmatrix} g\\r\end{pmatrix}$

\end{Prop}

\begin{pf}
We can give a coordinate description of the map 
$$\varphi_1: S^{g-1}H^*\rightarrow S^{g-2}H^*\otimes \wedge ^2H^*\otimes H$$
as follows. For any positive integer $k$, identify $S^kH$ with the
 space of polynomials of degree at most $k$ in a variable $x$
and $S^kH^*$ with the space of polynomials of degree at most $k$ in one
variable $y=x^*$. Identify $\wedge ^2H^*\rightarrow {\bf k}$ 
by the isomorphism $y\wedge 1\rightarrow 1$.
Then, the map  $\varphi_1$ acts by 
$$\varphi_1(y^k)=y^k\otimes x-y^{k-1} \otimes 1$$
with the convention that on the right hand side $y^{g-1}=0, y^{-1}=0$.
By taking the wedge products of this map, one then sees that
$$\varphi_r (y^{k_1}\wedge ...\wedge y^{k_r})=\sum_{0\le \epsilon_i\le 1}
(-1)^{\epsilon _1+...+\epsilon _r}y^{k_1-\epsilon_1}\wedge...
\wedge y^{k_r-\epsilon_r}
\otimes x^{r-(\epsilon _1+...+\epsilon _r)}$$
 again with the convention that on the right hand side $y^{g-1}=0,
 y^{-1}=0$.

A basis of $\wedge^r (S^{g-1}H^*)$ consists of the elements 
$y^{k_1}\wedge ...\wedge y^{k_r}$, $0\le k_1<...<k_r\le g-1$. The images
 of these elements are linearly independent as their leading terms (i.e.
the term of highest degrees on x and y jointly) obviously
are. Therefore, the dimension of $W^r$ is as stated.
\end{pf}

\section{General set up for infinitessimal deformations}
\begin{Rec}
\label{def}
We recall the basic set up for deformations of a curve, a vector bundle 
and its space of sections.
\end{Rec}
Write ${\bf k}[t]/t^2={\bf  k}_{\epsilon }$. By an infinitessimal
 defformation of the curve $C$ we mean a curve ${\cal C}_{\epsilon}$
over $Spec{\bf  k}_{\epsilon }$ with central fiber $C$. Similarly,
by an infinitessimal deformation of a vector bundle $F$ we mean
a vector bundle over $C\times Spec{\bf  k}_{\epsilon }$ with central
fiber $F$. By 
an infinitessimal deformation of the pair, we mean a curve
${\cal C}_{\epsilon}$ and a vector bundle ${\cal E}_{\epsilon}$
over ${\cal C}_{\epsilon}$.

Recall that the set of infinitessimal deformations of the curve $C$
can be parametrised by $H^1(C,T_C)$, the set of infinitessimal
deformations of the vector bundle $F$ can be parametrised by 
$H^1(F^*\otimes F)$ while the set of infinitessimal deformations of a pair
consisting of a curve $C$ and a vector bundle $F$ on $C$ can be parametrised 
by $H^1(\sum _F)$ where $\Sigma _F$ denotes the sheaf of first order 
differential operators acting on $F$.

We describe next the correspondence between these objects (cf.\cite{W}
proof of Prop. 1.2 and also \cite{BR} proof of 2.3). Assume given an element 
$\nu \in H^1(T_C)$. We think of the 
 sections of the sheaf  $T_C$ over an open set $U$ as the set of 
(${\bf k}$-linear)-maps ${\cal O}_U\rightarrow {\cal O}_U$ satisfying
$\nu (fg)=\nu (f)g+f\nu (g)$.
Take an affine open cover $C=\cup U_i$. Write $U_{ij}$ for $U_i\cap U_j$.
 Represent $\nu $ by  a 
cocycle $\nu =(\nu _{ij}), \nu _{ij}\in H^0(U_{ij},T_C)$. We associate 
to $\nu$ the following deformation of $C$: Consider the trivial 
deformations of the $U_i$, namely $U_i\times Spec {\bf k}_{\epsilon }$.
Glue them along the intersections $U_{ij}\times Spec {\bf k}_{\epsilon }$
using the matrices
$$\begin{pmatrix} Id&0\\
 \nu_{ij}&Id
\end{pmatrix}.$$
The correspondence 
$\nu _{ij}\rightarrow {\cal C}_{\epsilon}$
 obtained in this way is a bijection. 
 
Assume now given an element $\varphi \in H^1(F^*\otimes F)$.
 Represent it by a cocycle $(\varphi _{ij})$ with 
$\varphi _{ij} \in H^0(U_{ij},Hom(F,F))$. Consider the trivial extension 
of $F$ to $U_i\times Spec {\bf k}_{\epsilon }$, namely 
$F_{U_i}\oplus \epsilon F_{U_i}$. Take gluings on $U_{ij}$ given 
by 
$$\begin{pmatrix} Id&0\\
 \varphi _{ij}&Id
\end{pmatrix}.$$
This gives the correspondence between $H^1(F^*\otimes F)$
 and deformations of $F$.

 Assume now that a section $s$ of $F$ 
can be extended to a section $s_{\epsilon}$ of the deformation. 
There exist then local sections $s'_i \in H^0(U_i, F_{|U_i})$ 
such that $(s_{|U_i},s_i)$ define a section of $F_{\epsilon}$.
By construction of $F_{\epsilon}$ this means that
$$\begin{pmatrix} Id&0\\
 \varphi_{ij}&Id
\end{pmatrix}
\begin{pmatrix} s_{|U_i}\\
s'_i
\end{pmatrix}=
\begin{pmatrix} s_{|U_j}\\
s'_j
\end{pmatrix}.$$

This can be written as $\varphi_{ij}(s)=s'_j-s'_i$. 
Equivalently, 
$$\begin{matrix} \varphi _{ij}\in Ker (& H^1(F^*\otimes F)&
\rightarrow &H^1(F))\\
 &\nu _{ij}&\rightarrow &\nu_{ij}(s) \end{matrix}$$
This result can be formulated using the language of Brill-Noether Theory:
the set of infinitessimal deformations of the vector bundle $F$ that
 have sections deforming  a certain subspace $V\subset H^0(F)$
 consists of the orthogonal to the image  of the Petri map
$$P_V: V\otimes H^0(K\otimes F^*)\rightarrow H^0(K\otimes F\otimes F^*).
\leqno (\ref{def}.1)$$

Assume now given an element $\sigma \in H^1(\Sigma  _F)$.
We think of  $\Sigma _F(U)$ as the set
 of additive morphisms $\sigma :F(U)\rightarrow F(U)$ such that for a 
suitable element $\nu _{\sigma} \in T_C$, 
$\sigma (fs)= \nu_{\sigma } (f)s+f\sigma (s)$.
Represent $\sigma$ by a cocycle $\sigma =(\sigma _{ij}),\  \sigma_{ij}\in
H^1(U_{ij},\Sigma _F)$. Consider the associated element $(\nu _{ij})
\in H^1(T_C)$ and the corresponding deformation ${\cal C}_{\epsilon}$
of $C$. Take then the vector bundle on ${\cal C}_{\epsilon}$ obtained 
by gluing the trivial extensions of $F$ on $U_i$ by means of 
the matrices 
$$\begin{pmatrix} Id&0\\
 \sigma_{ij}&Id
\end{pmatrix}.$$

As in the case of deforming the line bundle alone, deformation of sections
is easy to interpret:
the set of infinitessimal deformations of the pair $(C,F)$ that have sections 
deforming  a certain subspace $V\subset H^0(F)$
consists of the orthogonal to the image  of the Petri map

$$\bar P_V:V\otimes H^0(K\otimes F^*)\rightarrow H^0(K\otimes \Sigma _F^*)
\leqno (\ref{def}.2)$$

defined as the dual of the natural cup-product map

$$H^1(\Sigma _F)\rightarrow Hom(V,H^1(F)).$$

Consider the exact sequence
$$
0\rightarrow F^*\otimes F\rightarrow \Sigma_F\rightarrow T_C\rightarrow 0.$$
The canonical map $\pi :\Sigma _F \rightarrow T_C$ is defined 
by $\pi (\sigma)=\nu _{\sigma}$. The map $i:F^*\otimes F\rightarrow \Sigma_F$
sends an element of $F^*\otimes F$ (considered as an endomorphism
of $F$) to itself. 
One obtains a commutative diagram
$$\begin{array}{cccccccccc}
 &0&\rightarrow &H^1(F^*\otimes F)&\rightarrow &H^1(\Sigma _F)
&\rightarrow&H^1(T_C)&\rightarrow &0\cr
(\ref{def}.3)& & &\downarrow P_V^*& &\downarrow \bar P_V^* & 
&\downarrow P_V^{'*} & & \cr
 &0&\rightarrow &ImP^*&\rightarrow &Hom(V,H^1(F))&\rightarrow&
(Ker P)^*&\rightarrow &0\cr 
\end{array}$$
and its dual (cf. \cite{AC}p.18 )
$$\begin{array}{cccccccccc}
 &0&\leftarrow&H^0(K\otimes F\otimes F^*)&\leftarrow
&H^0(K\otimes (\Sigma_F)^*)&\leftarrow &H^0(2K)&\leftarrow &0\cr
(\ref{def}.4)& & &\uparrow & &\uparrow \bar P_V& &\uparrow P'_V& & \cr
 &0&\leftarrow &ImP& \leftarrow &V \otimes H^0(K\otimes F^*))&
\leftarrow & Ker P&\leftarrow &0\cr 
\end{array}$$
When $V=H^0(F)$, we shall write $P_F$ instead of $P_V$. When $V$
and $F$ are clear, we shall suppress them from the notations.

We shall later use the following result. Its proof appears in \cite{T} 
Lemma 2.12.

\begin{Lem}
\label{1x-x1}
Let $M$ be a line bundle on a curve $C$ with two independent sections $s_0,s_1$
and such that $|K\otimes M^{-2}|$ has a section $t$. Denote by $D_1$
 the fixed part of the series determined by $s_0,s_1$,
 denote by $D_2$ the divisor corresponding to the section
$t$. Denote by $R$
the ramification divisor of the map $C\rightarrow {\bf P}^1$
 associated to the series $<s_0,s_1>$. 
Then, $P'_M(s_0\otimes ts_1-s_1\otimes ts_0)$ corresponds to the divisor
$2D_1+D_2+R$. In particular it is non-zero.
\end{Lem}

\begin{Lem}
\label{F=sum}
Assume that $F=\oplus _{i=1}^n F_i$ is a direct sum of vector bundles.
Then $$\Sigma _F=\oplus _{T_C}\Sigma_{F_i}\oplus [\oplus_{i\not= j}
F_i^*\otimes F_j]$$
Here $\oplus _{T_C}\Sigma_{F_i}$ denotes the fibered product over $T_C$
of the $\Sigma_{F_i}$.
\end{Lem}
\begin{pf}
 Consider an open set $U$. Let $\sigma :F(U)\rightarrow F(U)$ be 
a first order differential operator acting
on $F$. Using the decomposition of $F$ as a direct sum, $\sigma$ admits
 a representation as a matrix $(\sigma _{ij})$ where 
$\sigma _{ij}:F_i\rightarrow F_j$. Take a local section $s_k\in F_k(U)$.
One then checks that 
 $$\sigma ((0...0,fs_k,0...0))=(\sigma_{1k}(fs_k),...,\sigma_{nk}(fs_k))$$ 
Using that $\sigma ((0...0,fs_k,0...0))=f\sigma ((0...0,s_k,0...0))
+\nu (f)(0...0,s_k,0...0)$, we find that $\sigma _{kk}$ is a first order
differential operator corresponding to the same $\nu$ as $\sigma$
while $\sigma _{ij}$ is ${\cal O}_C$-linear if $i\not= j$.
\end{pf}

\begin{Lem}
\label{F=sum'}
Assume that $F=\oplus _{i=1}^n F_i$ is a direct sum of vector bundles.
 Denote by $P_i,\bar P_i,P'_i$ the Petri maps corresponding to the 
$F_i$. Then $P'_F$ can be obtained as the composition
$$Ker P\rightarrow \oplus _{i=1}^n KerP_i\rightarrow H^0(2K)$$
where the first map is the projection and the second map is ${1\over n}
\oplus P'_i$
 \end{Lem}
\begin{pf}
Consider the right hand square in \ref{def}.3 for each one of the $F_i$.
As in \ref{F=sum}, consider the fibered product of the $\Sigma_{F_i}$ over 
$T_C$. One then has a  commutative square
$$\begin{array}{ccc}
\oplus_{H^1(T_C)} H^1(\Sigma _{F_i})&\rightarrow&H^1(T_C)\cr
\downarrow \oplus_{i=1}^n\bar P_i^*& &\downarrow \oplus_{i=1}^n  P_i^{*} \cr
\oplus_{i=1}^nHom(H^0(F_i),H^1(F_i))&\rightarrow&\oplus_{i=1}^n coker P_i^{'*}
\cr 
\end{array}$$
 Take also the corresponding square for $F$ 
$$\begin{array}{ccc}
 H^1(\Sigma _F)&\rightarrow&H^1(T_C)\cr
\downarrow \bar P* & &\downarrow P^{'*} \cr
Hom(H^0(F),H^1(F))&\rightarrow&coker P^*\cr 
\end{array}$$
Consider the cube that has these diagrams as back and front faces 
respectively. We define four maps in the side edges.
The maps 
$$\oplus_{H^1(T_C)} H^1(\Sigma _{F_i})\rightarrow H^1(\Sigma _F) $$
and 
$$\oplus_{i=1}^nHom(H^0(F_i),H^1(F_i)) \rightarrow   Hom(H^0(F),H^1(F))$$
are natural diagonal injections (with zeroes on the terms corresponding
to a pair $F_i,F_j,i\not= j$). Notice that from \ref{F=sum}, the first map
is well defined. With these definitions, the left hand square commutes.

The map 
$$\oplus_{i=1}^n (Ker P_i)^* \rightarrow (Ker P)*$$
is defined as the dual of the natural projection. By dualisation, one can then 
check that the bottom face commutes.

If we take as the fourth map the homotethy
$$\times n :H^1(T_C) \rightarrow H^1(T_C)$$
 then, the top face commutes too.
 As $H^1(\Sigma _F)\rightarrow H^1(T_C)$ is onto,
 this shows that the right hand square commutes. Dualising this square, 
one obtains the result in the Lemma.

\end{pf}

\section{Deformations of E}

We apply the set up of the previous section to the
 hyperelliptic curve $C_0$ and 
the vector bundles $\wedge ^rE$.  As in section 1, $L$ denotes the
 hyperelliptic line bundle (i.e. the line bundle on $C$ of degree
two with two sections).
 We shall assume in all that follows 
that $r\le g-1-r$

\begin{Prop}
\label{P}
 The Petri map $P_{H^0(\wedge ^r E)}$ (cf. (\ref{def}.2))
associated to the vector bundle
 $\wedge ^r (E)$ gives by restriction an isomorphism 
$$P_{W^r}: W^r\otimes H^0(K\otimes (\wedge ^r E)^*)\rightarrow
 H^0(K\otimes \wedge^r E\otimes (\wedge ^r E)^*)$$
\end{Prop}
\begin{pf}
From \ref{W^r}
$$\wedge ^rE=\wedge ^r(S^{g-2}H^*\otimes \wedge ^2H^*)\otimes L^r$$
 so one has
 $$\wedge ^rE^*=\wedge ^r(S^{g-2}H\otimes \wedge ^2H)\otimes L^{-r}$$
Hence,
$$H^0(\wedge ^rE)=\wedge ^r(S^{g-2}H^*\otimes \wedge ^2H^*)\otimes 
S^rH.$$
and 
 $$H^0(K\otimes\wedge ^rE^*)=\wedge ^r(S^{g-2}H\otimes \wedge ^2H)\otimes 
S^{g-1-r}H.$$ 
Then, the Petri map $P_{H^0(E)}$ and its restriction $P_{W^r}$
to $W^r\otimes H^0(K\otimes (\wedge ^r E)^*) $ can be written as the tensor 
product with the vector space $ \wedge ^r(S^{g-2}H\otimes \wedge ^2H)$
of the diagram 
$$\begin{matrix}
\wedge ^r(S^{g-2}H^*\otimes \wedge ^2H^*)\otimes S^rH\otimes S^{g-1-r}H&
\rightarrow &\wedge ^r(S^{g-2}H^*\otimes \wedge ^2H^*)\otimes S^{g-1}H\cr
 \uparrow & &\uparrow \cr
W^r\otimes S^{g-1-r}H& \rightarrow &\wedge ^r(S^{g-2}H^*\otimes \wedge ^2H^*)
\otimes S^{g-1}H\cr \end{matrix}$$
Therefore, it is enough to show that the map $p_{W^r}$ in the lower row of
this diagram is an isomorphism.
Recall that $\varphi _r$ identifies $\wedge ^rS^{g-1}H^*$ to its image $W^r$
(where $\varphi _r$ is the natural immersion 
$\wedge ^rS^{g-1}H^*\rightarrow \wedge ^r(S^{g-2}H^*\otimes \wedge ^2H^*)
\otimes S^r$ defined in \ref{W^r}). The natural product map $S^rH\otimes 
S^{g-1-r}H\rightarrow S^{g-1}H$ can be identified to $P_{L^r}$.
Then, $p_{W^r}=(Id\otimes P_{L^r})o(\varphi_r\otimes Id)$.
We shall show that $p_{W^r}$ is an isomorphism by exhibiting its inverse
$q_{W^r}$.
Using the notations in \ref{W^r}, and the identification
$\wedge ^2H^*\cong {\bf k}$, one obtains
$$p_{W^r}(y^{k_1}\wedge ...\wedge y^{k_r}\otimes x^a)=
(Id\otimes P_{L^r})o(\varphi_r\otimes Id)
(y^{k_1}\wedge ...\wedge y^{k_r}\otimes x^a)=$$
$$=Id\otimes P_{L^r}
(\sum _{0\le \epsilon_i\le 1}(-1)^{\epsilon_1+...+\epsilon _r}
y^{k_1-\epsilon_1}\wedge ...\wedge 
y^{k_r-\epsilon_r}\otimes x^{r-(\epsilon _1+...+\epsilon _r)}\otimes x^a)=$$
$$=\sum _{0\le \epsilon_i\le 1}(-1)^{\epsilon_1+...+\epsilon _r}
y^{k_1-\epsilon_1}\wedge ...\wedge 
y^{k_r-\epsilon_r}\otimes x^{a+r-(\epsilon _1+...+\epsilon _r)}$$
with the convention that on the right hand side $y^{-1}=0,\ y^{g-1}=0$.
We describe $q_{W^r}$ as follows. Assume given integers  $0\le j_1<...<j_r
\le g-2,0\le b\le g-1$. There is then a value $l$ with $0\le l\le r$ such
that $j_l+1\le b\le j_{l+1}$. Define then 
$$q_{W^r}(y^{j_1}\wedge ...\wedge y^{j_r}\otimes x^b)=
\sum _{0\le t_i\le j_i-j_{i-1}-1,1\le s_i\le j_{i+1}-j_i}(-1)^{r-l}
y^{j_1-t_1}\wedge ...\wedge y^{j_l-t_l}\wedge y^{j_{l+1}+s_{l+1}}\wedge ...
y^{j_r+s_r}$$
$$\otimes x^{b-r-t_1-...-t_l+s_{l+1}+...+s_r}$$
with the conventions $j_{r+1}=g-1, j_{-1}=-1$.
Notice that the map $q_{W^r}$ is well defined as 
$0\le b-r-t_1-...-t_l+s_{l+1}+...+s_r\le g-1-r$ and 
$0\le j_1-t_1<...<j_l-t_l<j_{l+1}+s_{l+1}\le ...\le j_r+s_r\le g-1$.
It is a slightly tedious but straightforward computation to show that
the composition 
$p_{W^r}oq_{W^r}=Id$. As the two vector spaces involved have the same
 dimension, this suffices in order to prove the isomorphism.
\end{pf}

\begin{Cor}
\label{sigma_nu}
For any given infinitessimal deformation $\nu$ of the curve $C_0$, 
there is a unique infinitessimal deformation $\sigma$ of the 
pair consisting of $C_0$ and the vector bundle $\wedge ^rE$
such that $W^r$ can be extended to a space of sections of the deformation.
\end{Cor}
\begin{pf}

Consider the Petri map
 $$P_{W^r}:
W^r\otimes H^0(K\otimes F^*)\rightarrow H^0(K\otimes F\otimes F^*)$$
Consider the commutative diagram (\ref{def}.3) in case $F=\wedge ^rE, V=W^r$
$$\begin{array}{ccccccccc}
0&\rightarrow &H^1((\wedge ^rE)^*\otimes \wedge ^rE)&\rightarrow &H^1(\Sigma _
{\wedge ^rE})
&\rightarrow&H^1(T_C)&\rightarrow &0\cr
 & &\downarrow P^*& &\downarrow \bar P^* & &\downarrow P^{*'} & & \cr
0&\rightarrow &ImP^*&\rightarrow &Hom(W^r,H^1(\wedge ^r E))&\rightarrow&
(ker P)^*&\rightarrow &0\cr
\end{array}$$
From \ref{P},  $P^*$ is an isomorphism, $KerP=0$. Hence, every element in
 $H^1(T_C)$ has a 
unique inverse image in $H^1(\Sigma _F)$ that belongs to the kernel 
of $\bar P^*$. This proves the result.
\end {pf}

\begin{Cor}
\label{sigma_nu'}
For any given infinitessimal deformation $\nu$ of the curve $C_0$, the 
unique infinitessimal deformation $\sigma$ of the vector bundle
$\wedge ^rE$ that preserves it as the  exterior power of the dual of 
the kernel of the evaluation map is the $\sigma $ above.
\end{Cor}
\begin{pf} An infinitessimal deformation of $E$ as the dual of the kernel
of the evaluation map preserves the subspace $W=W^1$ as space of 
sections . Hence, an infinitessimal deformation of $\wedge ^rE$ as the 
exterior product of this dual preserves $W^r$ as space of sections. 
By the unicity of such deformation, the result follows.
\end{pf}

\begin{Prop}
\label{noW'}
Take $\nu$ a generic infinitessimal deformation of $C_0$. Let $\sigma$
be the deformation of $\wedge ^r E$ associated to $\nu$ as in \ref{sigma_nu'}.
If $\hat W$ is a subspace of $H^0(\wedge ^r(E))$ that strictly contains 
$W^r$, then $\hat W$ does not extend to a space of sections of 
the infinitessimal deformation
of $\wedge ^r E$ corresponding to $\sigma $.
\end{Prop}
\begin{pf}
It is enough to prove the result when $\hat W$ has dimension $a=dimW^r+1$.
Denote by $S$ the Schubert cycle
of subspaces of dimension $a$ of $H^0(\wedge ^rE)$ that contain 
$W^r$ ($S\subset {\bf Gr}(a, H^0(\wedge ^rE))$)  .

Consider the diagram (\ref{def}.4) for the case $F=\wedge ^rE, V=\hat W$.
 We obtain
$$\begin{array}{ccccccccc}
0&\leftarrow&H^0(K\otimes (\wedge ^rE)\otimes (\wedge ^rE)^*)&\leftarrow
&H^0(K\otimes (\Sigma_{\wedge^rE})^*)&\leftarrow &H^0(2K)&\leftarrow &0\cr
 & &\uparrow & &\uparrow \bar P_{\hat W}& &\uparrow P'_{\hat W}& & \cr
0&\leftarrow &ImP& \leftarrow &\hat W\otimes H^0(K\otimes (\wedge ^rE)^*)&
\leftarrow & Ker P_{\hat W}&\leftarrow &0\cr
\end{array}$$
Define $P'=P'_{H^0(\wedge ^r(E))}$.

We first show that the result would follow from the following statement:

(*) There is an element $\alpha \in H^0(2K)$ such that $\alpha \in P'(\hat W)$
for all $\hat W in S$.

If (*) holds, consider a direction of deformation $\nu$
in $H^1(T_C)$ not orthogonal to $\alpha$ and its corresponding
deformation $\sigma_{\nu}$. Then, $\bar P^*_{\hat W}(\sigma _{\nu })\not= 0$
 for 
all $\hat W$. Therefore, no $\hat W$  extends to a space of sections 
of the deformation corresponding to $\sigma _{\nu}$.

Let us see that in fact it suffices to prove (*) for the generic 
 $\hat W \in S$.
Note then that from \ref{P}, $Ker P_{\hat W}$ is a subspace of
 $H^0(\wedge ^rE)\otimes H^0(K\otimes (\wedge ^rE)^*)$ of fixed dimension
 $b=h^0(K\otimes (\wedge ^rE)^*)$.  We have a well-defined map

$$\begin{array}{cccc} \varphi : &S&\rightarrow &
{\bf Gr}(b,H^0(\wedge ^rE)\otimes H^0(K\otimes (\wedge ^rE)^*)\cr
 &\hat W&\rightarrow &KerP_{\hat W}
\end{array}$$

Define $\bar P=\bar P_{H^0(\wedge ^r(E))}$. Assume that for the generic 
$\hat W\in S$, $\varphi (\hat W)$ intersects
$(P')^{-1}(\alpha)$ for a fixed element $\alpha \in H^0(2K)$.
Denote by $S'$ the Schubert cycle of $ {\bf Gr}(b,H^0(\wedge ^rE)
\otimes H^0(K\otimes (\wedge ^rE)^*)$ of subspaces that intersect
$(P')^{-1}(\alpha)$. Our assumption is that for $\hat W\in S$
generic, $\varphi (\hat W)\in S'$. From the irreducibility of $S$
and  the closedness of $S'$, the same is true for every element in  
$S$.

 Hence, it remains to prove (*) for generic $\hat W\in S$.
This  will follow from the Lemma below and \ref{1x-x1}.
\end{pf}
\begin{Lem}
\label{1x-x1'}
Let $\hat W$ be generic in $S$, then
 $P'_{\hat W}(Ker P_{\hat W})=P'_{L^r}(Ker P_{L^r})  $.
\end{Lem}

\begin{pf}[of \ref{1x-x1'}]
As noticed before, $dimKer P_{\hat W}$ is constant for $\hat W\in S$.
 Therefore, the map
$P'_{\hat W}$ attains its maximum rank on an open set of $S$.
From \ref{F=sum'}, the image of $P'_{\hat W}$ is contained in the 
image of $P'_{L^r}$. Hence, the maximum possible rank of $P'_{\hat W}$
is $dim ImP'_{L^r}$. It then suffices 
to exhibit one $\hat W$ satisfying the condition.
 We use the basis of $\wedge ^r (S^{g-2}H^*\otimes \wedge ^2H^*)$ introduced
 in \ref{W^r} and the corresponding trivialisation of $\wedge ^2H^*$.
Define 
$$\hat w=\sum _{j=0}^{r-1}(-1)^j1\wedge y\wedge ...\wedge y^{j-1}\wedge y^{j+1}
\wedge ...\wedge y^r\otimes x^{r-j}$$
Let $\hat W$ be the space generated by $\hat w$ and $W^r$.
We recall from the proof of \ref{P} that the map $P_{\hat W}$ can be
 obtained by taking tensor product with 
$ \wedge ^r(S^{g-2}H\otimes \wedge ^2H)$
of the map 
$$p_{\hat W}:\hat W\otimes S^{g-1-r}H\rightarrow 
\wedge ^r(S^{g-2}H^*\otimes \wedge ^2H^*)\otimes S^{g-1}H$$
 Consider the 
following element of $\hat W\otimes S^{g-1-r}H$  
$$z=\hat w\otimes 1+\sum _{i=r+1}^{g-1}\varphi_r(y\wedge y^2\wedge ...
\wedge y^{r-1}\wedge y^i)\otimes x^{i-r}$$
We want to see that $p_{\hat W}(z)=0$.

We can write an explicit expression for $z$ using the convention
that $y^{-1}=0,y^{g-1}=0$
$$z=\sum _{j=0}^{r-1}(-1)^j1\wedge y\wedge ...\wedge y^{j-1}\wedge y^{j+1}
\wedge ...\wedge y^r\otimes x^{r-j}\otimes 1$$
$$+\sum_{i=r+1}^{g-1}\sum_{0\le \epsilon_i\le 1}
(-1)^{\epsilon _1+...+\epsilon _r}y^{1-\epsilon_1}\wedge...
\wedge y^{r-1-\epsilon_{r-1}}\wedge y^{i-\epsilon _r}\otimes 
 x^{r-(\epsilon _1+...+\epsilon _r)}\otimes x^{i-r}$$
$$=\sum _{j=0}^{r-1}(-1)^j1\wedge y\wedge ...\wedge y^{j-1}\wedge y^{j+1}
\wedge ...\wedge y^r\otimes [x^{r-j}\otimes 1-x^{r-(j+1)}\otimes x]$$
$$+\sum_{i= r+1}^{g-2}\sum_{0\le \epsilon_i\le 1}
(-1)^{\epsilon _1+...+\epsilon _{r-1}}y^{1-\epsilon_1}\wedge...
\wedge y^{r-1-\epsilon_{r-1}}\wedge y^i\otimes$$ 
$$[ x^{r-(\epsilon _1+...+\epsilon _{r-1})}\otimes x^{i-r}-
 x^{r-(\epsilon _1+...+\epsilon _{r-1}+1)}\otimes x^{i+1-r}]$$
It is clear from the above expression for $z$ that $z\in Kerp_{W'}$.
In particular, for any element $e\in \wedge ^r(S^{g-2}H\otimes \wedge ^2H)$,
$z\otimes e \in KerP_{W'}$.  

From \ref{F=sum'}, the map $P'_{\hat W}$ is up to a constant $c'$
($c'={1\over { g-2\choose r}}$) the composition of the
following maps: 
natural inclusion 
$$\hat W\otimes \wedge ^r(S^{g-2}H\otimes \wedge ^2H)\otimes S^{g-1-r}H
  \rightarrow $$
$$\wedge ^r(S^{g-2}H^*\otimes \wedge ^2H^*)\otimes
\wedge ^r(S^{g-2}H\otimes \wedge ^2H)\otimes S^{g-1}H    $$
 the cup-product with the identity element
in $\wedge ^r(S^{g-2}H^*\otimes \wedge ^2H^*) \otimes
\wedge ^r(S^{g-2}H\otimes \wedge ^2H)$
$$\wedge ^r(S^{g-2}H^*\otimes \wedge ^2H^*)\otimes
\wedge ^r(S^{g-2}H\otimes \wedge ^2H \otimes S^{g-1}H  
\rightarrow S^rH\otimes S^{g-1-r}H$$
followed by the natural product map 
$$P'_{L^r}:S^rH\otimes S^{g-1-r}H\rightarrow S^{g-1}H$$.

 Hence, for $j=0,...,r-1$ writting $c=(-1)^jc'$,
$$P'_{\hat W}(z\otimes 1\wedge x\wedge ...\wedge x^{j-1}\wedge x^{j+1}
\wedge ...\wedge x^r)=cP'_{L^r}[x^{r-j}\otimes 1-x^{r-(j+1)}\otimes x]$$
and for $i=r+1...g-2, 0\le \epsilon _i\le 1$ writting
 $c=(-1)^{\epsilon_1+...+\epsilon_{r-1}}c'$
$$P'_{\hat W}(z\otimes x^{1-\epsilon_1}\wedge...\wedge x^{r-1-\epsilon_{r-1}}
\wedge x^i )=
cP'_{L^r}[x^{r-(\epsilon _1+...+\epsilon _{r-1})}\otimes x^{i-r}-
 x^{r-(\epsilon _1+...+\epsilon _{r-1}+1)}\otimes x^{i+1-r}]$$
As $ker P_{L^r}$ is generated by elements of the form 
$$x^i\otimes x^j-x^{i+1}
\otimes x^{j-1}, \  i=0...r-1,\ j=1...g-1-r$$
the result is proved.
\end{pf}

\begin{section}{Extending the results to the generic curve.}

\begin{Prop} 
\label{inj}
Let $C$ be a generic curve of genus $g$. Let $E$ be defined as in \ref{E}.
Denote by $W$ the image of $(H^0(C,K_C))^*$ in $H^0(C,E)$. Then,
 the natural map 
$$\psi _{C,r}: \wedge ^rW \rightarrow H^0(C,\wedge ^rE)$$ 
is injective.
\end{Prop}
Note that for $C$ non-hyperelliptic, $W=H^0(E)$.
This follows from the projective normality of $C$ (case $p=0$ of the 
conjecture).

\begin{pf}
If for a given curve $C$, $\psi _{C,r}$ is injective, the same holds 
for every curve in a neighborhood of $C$ in ${\cal M}_g$.
As $dim \wedge ^rW=\binom{g}{r}$, \ref{W^r} shows that $\psi _{C_0,r}$
is injective for $C_0$ hyperelliptic. Hence the result follows.
\end{pf}

The following proposition now concludes the proof of \ref{teorema}:

\begin{Prop}
Let $C$ be a generic curve of genus $g$. Then, $h^0(C,\wedge ^rE)=\binom{g}
{r}$ and $\wedge ^rW\rightarrow H^0(C,\wedge ^rE)$ is an isomorphism.
\end{Prop}
\begin{pf}
From \ref{inj}, $h^0(C,\wedge ^rE)\ge dim Im
(\wedge ^rW \rightarrow H^0(C,\wedge ^rE))=dim\wedge ^rW=\binom{g}{r}$.
From \ref{noW'}, $h^0(C,\wedge^rE)\le dim W^r_{C_0}=\binom{g}{r}$.
This concludes the proof.
\end{pf}
\end{section}

\end{document}